\documentclass[11pt,letterpaper,reqno]{amsart}
\usepackage{amsmath,amsthm,amsfonts,amssymb,mathtools,algpseudocode}
\usepackage{comment,bm}

\usepackage{bookmark,hyperref}
\hypersetup{pdfstartview={FitH}}

\newtheorem{Theorem}{{\bf Theorem}}

\newtheorem{Corollary}[Theorem]{{\bf Corollary}}

\newtheorem{Proposition}[Theorem]{{\bf Proposition}}
\newtheorem{Definition}[Theorem]{{\bf Definition}}

\newtheorem{Remark}[Theorem]{{\bf Remark}}

\newcommand{\Ol}{\mbox{\Large $\mathcal O$}}

\newcommand{\jd}{\mathsf{12}}
\newcommand{\jt}{\mathsf{13}}
\newcommand{\jc}{\mathsf{14}}
\newcommand{\dt}{\mathsf{23}}
\newcommand{\dc}{\mathsf{24}}
\newcommand{\tc}{\mathsf{34}}

\def\sh{\mathop{\rm sh}\nolimits}
\def\ch{\mathop{\rm ch}\nolimits}
\def\th{\mathop{\rm th}\nolimits}

\numberwithin{equation}{section}

\begin{document}

\title[Jacobi method for symmetric $4\times4$ matrices]{Jacobi method for symmetric $4\times4$ matrices converges for every cyclic pivot strategy}

\author{Erna Begovi\'{c}~Kova\v{c}}\thanks{Erna Begovi\'{c} Kova\v{c}, Faculty of Chemical Engineering and Technology, University of Zagreb, Maruli\'{c}ev trg 19, 10000 Zagreb, Croatia, \texttt{ebegovic@fkit.hr}}
\author{Vjeran Hari}\thanks{Vjeran Hari, Department of Mathematics, University of Zagreb, Bijeni\v{c}ka 30, 10000 Zagreb, Croatia, \texttt{hari@math.hr}}
\thanks{This work has been fully supported by Croatian Science Foundation under the project 3670.}
\date{1 August 2017}


\begin{abstract}
The paper studies the global convergence of the Jacobi method for symmetric matrices of size $4$. We prove global convergence for all $720$ cyclic pivot strategies.
Precisely, we show that inequality $S(A^{[t+3]})\leq\gamma S(A^{[t]})$, $t\geq1$, holds
with the constant $\gamma<1$ that depends neither on the matrix $A$ nor on the pivot strategy.
Here $A^{[t]}$ stands for the matrix obtained from $A$ after $t$ full cycles of the Jacobi method and $S(A)$ is the off-diagonal norm of $A$.
We show why three consecutive cycles have to be considered.
The result has a direct application on the $J$-Jacobi method.
\end{abstract}

\maketitle

\section{Introduction}

This paper considers the global convergence of the Jacobi method and is a continuation of the work from~\cite{BH_4x4pa2015,BH_Jac2015}. It answers the question whether all cyclic Jacobi methods are convergent, at least for the case $n=4$. Furthermore, the result  implies that any cyclic $J$-Jacobi method (see~\cite{Ves-JJac93,HSS14}) for solving the generalized eigenvalue problem $Ax=\lambda Jx$ is globally convergent for $n=4$.
The global convergence of the cyclic Jacobi methods for symmetric matrices has been considered in~\cite{for+hen-60,han-63,hen+zim-68,N75,SS89,mas-95} and lately in~\cite{B_PhD,BH_4x4pa2015,H15,BH_Jac2015}. These latest results have sufficiently enlarged the class of ``convergent'' strategies to make the general proof in the case $n=4$ possible.

Our results have a direct application on the corresponding block Jacobi methods with $2\times2$ blocks. In that case pivot submatrices are $4\times4$ matrices. In order to have a convergent block method we need a convergent core algorithm. Common approach in each block step of the block Jacobi method is to use the spectral decomposition of the pivot submatrix. This paper shows that, in the case of the blocks of size $2$, one can use the Jacobi algorithm with an arbitrary cyclic pivot strategy as the core algorithm. This choice of the core algorithm is natural because of the quadratic convergence and high relative accuracy of the Jacobi algorithm. In addition, it ensures the global convergence of the block Jacobi method (see~\cite{BH_Jac2015}, Theroem~2.10(ii)).

The eigenvalue problems for symmetric and $J$-symmetric matrices of size $4$ can be solved directly, by forming the quartic characteristic polynomial and using Cardano's rootfinding formulas. However, efficiency and stability of such approach can be questioned, since there are fast and stable algorithms based on Jacobi and $J$-Jacobi methods. Especially in
the context of block Jacobi methods with $2\times2$ blocks when the pivot submatrices are of size $4$ and most of the time they are nearly diagonal. On such matrices Jacobi methods are very efficient (just few sweeps are needed) and highly accurate \cite{mat-09}.

Jacobi method for symmetric matrices is an iterative process of the form
\begin{equation}\label{Jacobi}
A^{(k+1)}=R_k^TA^{(k)}R_k, \quad k\geq0, \qquad A^{(0)}=A=A^T,
\end{equation}
where $R_k=R(i(k),j(k),\phi_k)$, $k\geq 0$, are plane rotations. Here $R(p,q,\psi)$,
$p<q$, stands for the rotation in the $(p,q)$ plane, whose $(p,q)$-element is equal to $-\sin\psi$.
The goal of the $k$-th step is to reduce the off-diagonal norm of the current matrix $A^{(k)}$. The off-diagonal norm, or shortly off-norm, of a symmetric $n\times n$ matrix $ X=(x_{rs})$ is defined by
$$S(X)=\frac{\sqrt{2}}{2}\|X-\text{diag}(X)\|_F=\sqrt{\sum_{i=1}^{n-1}\sum_{j=i+1}^n x_{ij}^2},$$
where $\|X\|_F$ denotes the Frobenius norm of $X$. The rotation angle $\phi_k$ is defined by
$$\tan(2\phi_k)=\frac{2a_{ij}^{(k)}}{a_{ii}^{(k)}-a_{jj}^{(k)}}, \quad -\frac{\pi}{4}\leq\phi_k\leq\frac{\pi}{4}, \ k\geq0,$$
and this choice results in the annihilation of the \emph{pivot element} $a_{ij}^{(k)}$ of
$A^{(k)}$. Throughout this paper it is presumed that the Jacobi method is defined by the rotation angles from the interval $[-\frac{\pi}{4},\frac{\pi}{4}]$.
The pair $(i,j)$, $i<j$, depending on $k$, is called \emph{pivot pair}. Hence $i=i(k)$, $j=j(k)$, and letters $i$ and $j$ will be reserved for pivot indices.
A method is said to be \emph{globally convergent} if, for any initial $A$, the sequence of matrices $(A^{(k)},k\geq0)$ generated by the iterative process~\eqref{Jacobi}
converges to some diagonal matrix.

The diagonal elements of the iteration $A^{(k)}$ generated by the Jacobi method always converge~\cite{mas-95}. Hence, a method is globally convergent if and only if
the sequence of the off-norms $(S(A^{(k)}),k\geq0)$ converges to zero for any initial $A$. This sequence is nonincreasing because
$$S^2(A^{(k+1)})=S^2(A^{(k)})-(a_{ij}^{(k)})^2, \quad k\geq0.$$
Therefore, it is sufficient to find a subsequence of $(S(A^{(k)}),k\geq0)$ that converges to zero.
In particular, it is sufficient to prove that
\begin{equation}\label{intro_ineq}
S(A^{((t+\tau)N)})\leq\gamma S(A^{(tN)}), \quad t\geq t_0, \quad N=\frac{n(n-1)n}{2}
\end{equation}
holds for every cyclic strategy and every initial $A$, with $\tau,t_0\in\mathbb{N}$ and the constant $\gamma<1$, that depends neither on $A$ nor on the pivot strategy.
In the case $n=4$, we shall prove that relation~\eqref{intro_ineq} holds with $\tau=3$, $t_0=1$ and $\gamma=1-10^{-5}$.

For a symmetric $4\times4$ matrix there are altogether $720$ cyclic pivot strategies. The convergence properties of the Jacobi method greatly vary for different pivot strategies \cite{BH_4x4pa2015}. Say, for some parallel strategies from Section~3.3 and some initial symmetric matrices, the quotient $S(A^{(N)})/S(A)$ can be arbitrary close to $1$
\cite[Proposition~5.2]{BH_4x4pa2015}, and in the second cycle that quotient can be as small as $10^{-26}$ \cite[Example~5.1]{BH_4x4pa2015}. Hence, for convergence considerations one has to  take into account more than one cycle.

The paper is divided into $5$ sections.
In Section~\ref{sec:strategies} we present a brief theory of pivot strategies and pivot orderings. The main result of the paper is formulated in Theorem~\ref{tm:main}.
In Section~\ref{sec:strategies3} different classes of cyclic pivot strategies are analyzed. Specifically, the generalized serial pivot strategies from~\cite{B_PhD,BH_Jac2015} are studied in Subsections~\ref{sec:serialperm} and~\ref{sec:generalizedperm} and the parallel strategies from~\cite{B_PhD,BH_4x4pa2015} in Subsection~\ref{sec:parallel}.
In Section~\ref{sec:proof} the proof of the main theorem is presented. Finally, in Section~\ref{sec:appl} we show how the new results can be applied to the $J$-Jacobi method.

Most of the results presented here are a part of the unpublished thesis \cite{B_PhD}.

\section{Cyclic pivot strategies and pivot orderings}\label{sec:strategies}

In the $k$-th step of the Jacobi method, the pivot pair $(i(k),j(k))$ is selected and the pivot element $a_{ij}^{(k)}$ of $A^{(k)}$ is annihilated. The rule for choosing the pivot pairs is called the \emph{pivot strategy}. Pivot strategy can be identified with a function
$I:\mathbb{N}_0\rightarrow\mathbf{P}_n$,
where $\mathbb{N}_0=\{0,1,2,3,\ldots\}$, $\mathbf{P}_n=\big{\{}(r, s) \ \big{|} \ 1\leq r<s\leq n\big{\}}$, and $n$ is the size of the matrix.
If there is $T\in\mathbb{N}$ such that $I(k+T)=I(k)$, $k\geq0$,
then $I$ is \emph{periodic} with period $T$. In addition, if  $T=N\equiv\frac{n(n-1)}{2}$ and $\{I(k) \ \big{|} \ 0\leq k\leq T-1\}=\mathbf{P}_n$, the pivot strategy is \emph{cyclic}. The pivot pair $I(k)$ will be briefly denoted by $(i,j)$ when $k$ is implied.

Let $\Ol(\mathbf{P}_n)$ denote the set of all orderings of $\mathbf{P}_n$.
For a given pivot strategy $I$, the corresponding pivot ordering $\mathcal{O}_I\in\Ol(\mathbf{P}_n)$ is defined as
\begin{equation}\label{ordering}
\mathcal{O}_I=I(0),I(1),\ldots,I(N-1).
\end{equation}
It corresponds to one \emph{cycle} or \emph{sweep} of the method.
On the other hand, for a given ordering $\mathcal{O}=(i_0,j_0),(i_1,j_1),\ldots ,(i_{N-1},j_{N-1})$ from $\mathcal{\Ol}(\mathbf{P}_n)$, the corresponding pivot strategy $I_{\mathcal{O}}$ is defined by
$$I_{\mathcal{O}}(k)=(i_{\tau(k)},j_{\tau(k)}), \quad 0\leq\tau(k)\leq N-1, \ k\equiv\tau(k)(\bmod \ N), \ k\geq0.$$
Hence, $\mathcal{O}_{I_{\mathcal{O}}}=\mathcal{O}$ and we can use the terms ``ordering'' and ``pivot ordering'' as synonyms.

We say that the ordering $\mathcal{O}\in\mathcal{\Ol}(\mathbf{P}_n)$ and the pivot strategy $I_{\mathcal{O}}$ are \emph{convergent} if relation~\eqref{intro_ineq} holds for $I_{\mathcal{O}}$.

An \emph{admissible transposition} on $\mathcal{O}\in\mathcal{\Ol}(\mathbf{P}_n)$ is a transposition of two adjacent terms
$(i_r,j_r),$ $(i_{r+1},j_{r+1})\rightarrow(i_{r+1},j_{r+1}),(i_r,j_r)$, $0\leq r<N-1,$
provided that $\{i_r,j_r\}$ and $\{i_{r+1},j_{r+1}\}$ are disjoint. We say that such pairs $(i_r,j_r)$ and $(i_{r+1},j_{r+1})$ \emph{commute}.

Two orderings $\mathcal{O},\mathcal{O}'\in\mathcal{\Ol}(\mathbf{P}_n)$ are:
\begin{itemize}
\item[(i)] \emph{equivalent} (we write $\mathcal{O}\sim\mathcal{O}'$) if one can be obtained from the other using a finite number of admissible transpositions.
\item[(ii)] \emph{shift-equivalent} ($\mathcal{O}\stackrel{s}{\sim}\mathcal{O}'$) if $\mathcal{O}=[\mathcal{O}_1,\mathcal{O}_2]$ and $\mathcal{O}'=[\mathcal{O}_2,\mathcal{O}_1]$, where $[ \ , \  ]$ stands for the concatenation. The number of pairs in the sequence $\mathcal{O}_1$ is called the length of the shift.
\item[(iii)] \emph{weakly equivalent} ($\mathcal{O}\stackrel{w}{\sim}\mathcal{O}'$) if there are $\mathcal{O}_i\in\mathcal{\Ol}(\mathbf{P}_n)$, $1\leq i\leq r$, such that
    in the sequence $\mathcal{O},\mathcal{O}_1,\ldots,\mathcal{O}_r,\mathcal{O}'$ every two adjacent terms are equivalent or shift-equivalent.
\end{itemize}
It is easy to check that $\sim$, $\stackrel{s}{\sim}$ and $\stackrel{w}{\sim}$ are equivalence relations on the set $\mathcal{\Ol}(\mathbf{P}_n)$.
Two pivot strategies $I_{\mathcal{O}}$ and $I_{\mathcal{O}'}$ are equivalent (shift-equivalent, weakly equivalent) if the related pivot orderings from relation~\eqref{ordering} are equivalent (shift-equivalent, weakly equivalent).

We will need two more relations between pivot orderings. They have been introduced in \cite{B_PhD,BH_Jac2015}. These relations will be used to prove the global convergence of the Jacobi method under different pivot orderings.

\begin{Definition}
Let $\mathcal{O}\in\mathcal{\Ol}(\mathbf{P}_n)$, $\mathcal{O}=(i_0,j_0),(i_1,j_1),\ldots,(i_{N-1},j_{N-1})$. Then
$$\mathcal{O}^{\leftarrow}=(i_{N-1},j_{N-1}),\ldots,(i_1,j_1),(i_0,j_0) \ \in\mathcal{\Ol}(\mathbf{P}_n)$$
is the reverse ordering to $\mathcal{O}$. The pivot strategy $I_{\mathcal{O}}^{\leftarrow}=I_{\mathcal{O}^{\leftarrow}}$ is reverse to $I_{\mathcal{O}}$.
\end{Definition}

\begin{Definition}
Let $\mathcal{O}\in\mathcal{\Ol}(\mathbf{P}_n)$, $\mathcal{O}=(i_0,j_0),(i_1,j_1),\ldots,(i_{N-1},j_{N-1})$. Then $\mathcal{O}'\in\mathcal{\Ol}(\mathbf{P}_n)$ is
\emph{permutation equivalent} to $\mathcal{O}$ if there is a permutation $\mathsf{q}$ of the set $\{1,2,\ldots,n\}$ such that
\begin{align*}
\mathcal{O}' & = \big{(}\min\{\mathsf{q}(i_0),\mathsf{q}(j_0)\},\max\{\mathsf{q}(i_0),\mathsf{q}(j_0)\}\big{)},\ldots \\
& \quad \ldots,\big{(}\min\{\mathsf{q}(i_{N-1}),\mathsf{q}(j_{N-1})\},\max\{\mathsf{q}(i_{N-1}),\mathsf{q}(j_{N-1})\}\big{)}.
\end{align*}
We write $\mathcal{O}'=\mathcal{O}(\mathsf{q})$ or $\mathcal{O}\stackrel{\mathsf{q}}{\sim}\mathcal{O}'$. The general notation $\mathcal{O}'\stackrel{p}{\sim}\mathcal{O}$ will
stand for permutation equivalence between $\mathcal{O}'$ and $\mathcal{O}$ with some permutation.
\end{Definition}

Let us now consider the case $n=4$. We have
$$\mathbf{P}_4=\big{\{}(1,2),(1,3),(1,4),(2,3),(2,4),(3,4)\big{\}} \quad \text{and} \quad N=6.$$
Since the cardinality of the set $\mathbf{P}_4$ is  $6$, there are $6!$ orderings of $\mathbf{P}_n$ and consequently $720$ cyclic pivot strategies. The ones of our special interest are those (cyclic) strategies that start with the pivot element at position $(1,2)$, i.e.\ for which $I(0)=(1,2)$. Set
$$\mathcal{C}_0=\mathcal{C}_0(\mathbf{P}_4):=\big{\{}\mathcal{O}\in\Ol(\mathbf{P}_4) \ \big{|} \ I_{\mathcal{O}}(0)=(1,2)\big{\}}.$$
There are $5!=120$ pivot orderings in $\mathcal{C}_0$. The global convergence results for all strategies $I_{\mathcal{O}}$, $\mathcal{O}\in\Ol(\mathbf{P}_4)$, will follow from the results obtained for $I_{\mathcal{O}}$, $\mathcal{O}\in\mathcal{C}_0$.

The main result of the paper is given in the following theorem.

\begin{Theorem}\label{tm:main}
Let $A$ be a symmetric $4\times4$ matrix and let $\mathcal{O}\in\Ol(\mathbf{P}_4)$. Let $A^{[t]}$ be obtained from $A$ by applying $t$ sweeps of the cyclic Jacobi method defined by the pivot strategy $I_{\mathcal{O}}$. Then there is a constant $\gamma$, $0\leq\gamma<1$, that depends neither on $A$ nor on the pivot strategy, such that
$$S(A^{[t+3]})\leq\gamma S(A^{[t]}), \quad t\geq 1.$$
\end{Theorem}

\begin{proof}
The proof uses results from Section~\ref{sec:strategies3}, and is moved to Section~\ref{sec:proof}.
\end{proof}

\section{Auxiliary results}\label{sec:strategies3}

Special classes of cyclic pivot strategies have been introduced and studied in~\cite{B_PhD,BH_Jac2015}. We briefly recall notation and convergence results from~\cite{BH_Jac2015}. The results refer to the Jacobi method for a symmetric $n\times n$ matrix. Later on, we will consider the case $n=4$.

\subsection{Classes of pivot orderings $\mathcal{C}_c^{(4)}$ and $\mathcal{C}_r^{(4)}$}\label{sec:serialperm}

Let us recall definitions of the sets $\mathcal{C}_c^{(n)}$ and $\mathcal{C}_r^{(n)}$ from~\cite{BH_Jac2015}. Let $\Pi^{(l_1,l_2)}$ denote the set of all permutations of the set $\{l_1,l_1+1,\ldots\ ,l_2\}$, $l_1\leq l_2$.
Then,
\begin{align}
\mathcal{C}_c^{(n)} = \big{\{}\mathcal{O}\in\mathcal{\Ol}(\mathbf{P}_n) \ \big{|} \ & \mathcal{O}= (1,2),(\pi_{3}(1), 3),(\pi_{3}(2), 3),\ldots,(\pi_{n}(1),n),\ldots,(\pi_{n}(n-1),n), \nonumber \\
& \pi_{j}\in\Pi^{(1,j-1)}, \ 3\leq j\leq n\big{\}} \label{Ccdef}
\end{align}
is the set of \emph{column-wise orderings with permutations} of $\mathcal{\Ol}(\mathbf{P}_n)$. The columns are taken from left to right.
Similarly,
\begin{align}
\mathcal{C}_r^{(n)} = \big{\{}\mathcal{O}\in\mathcal{\Ol}(\mathbf{P}_n) \ \big{|} \ & \mathcal{O}= (n-1,n),(n-2,\tau_{n-2}(n-1)),(n-2,\tau_{n-2}(n)),\ldots,(1,\tau_{1}(n)), \nonumber\\
& \tau_{i}\in\Pi^{(i+1,n)}, \ 1\leq i\leq n-2\big{\}} \label{Crdef}
\end{align}
is the set of \emph{row-wise orderings with permutations} of $\mathcal{\Ol}(\mathbf{P}_n)$. Here, the rows are taken from bottom to top.

The sets
\begin{equation*}\label{Cinvdef}
\overleftarrow{\mathcal{C}_c}^{(n)} = \big{\{} \mathcal{O}\in\mathcal{\Ol}(\mathbf{P}_n) \ \big{|} \ \mathcal{O}^{\leftarrow}\in \mathcal{C}_c^{(n)} \big{\}}, \quad
\overleftarrow{\mathcal{C}_r}^{(n)} = \big{\{} \mathcal{O}\in\mathcal{\Ol}(\mathbf{P}_n) \ \big{|} \ \mathcal{O}^{\leftarrow}\in \mathcal{C}_r^{(n)} \big{\}},
\end{equation*}
contain the orderings reverse to those from $\overleftarrow{\mathcal{C}_c}^{(n)}$, $\overleftarrow{\mathcal{C}_r}^{(n)}$, respectively. This means that in the orderings from $\overleftarrow{\mathcal{C}_c}^{(n)}$ ($\overleftarrow{\mathcal{C}_r}^{(n)}$) the pairs are arranged column-by-column from right to left (row-by-row from top to bottom), with permutations inside the columns (rows).

The union of all these sets is denoted by $\mathcal{C}_{sp}^{(n)}$,
$$\mathcal{C}_{sp}^{(n)}=\mathcal{C}_c^{(n)}\cup\overleftarrow{\mathcal{C}_c}^{(n)}
\cup\mathcal{C}_r^{(n)}\cup\overleftarrow{\mathcal{C}_r}^{(n)},$$
and is called the set of \emph{serial orderings with permutations} of $\mathcal{\Ol}(\mathbf{P}_n)$.

The convergence result for the Jacobi method defined by the pivot strategies $I_{\mathcal{O}}$, $\mathcal{O}\in\mathcal{C}_{sp}^{(4)}$, is given in Proposition~\ref{tm:sp}. This is a special case of~\cite[Corollary 3.3]{BH_Jac2015} for $n=4$.

\begin{Proposition}\label{tm:sp}
Let $A$ be a symmetric $4\times4$ matrix and $\mathcal{O}\in\mathcal{C}_{sp}^{(4)}$. Let $A'$ be obtained from $A$ by applying one sweep of the cyclic Jacobi method defined by the strategy $I_{\mathcal{O}}$. Then
\begin{equation}\label{sp}
S^2(A')\leq \eta_4 S^2(A), \quad \eta_4=\frac{27}{28}.
\end{equation}
\end{Proposition}

\begin{proof}
It has been shown in~\cite[Corollary 3.3]{BH_Jac2015} that the relation (\ref{sp}) holds for a symmetric matrix $A$ of size $n$ and $\mathcal{O}\in\mathcal{C}_{sp}^{(n)}$ with the
constant $\eta_n$ satisfying
\begin{equation}\label{3.5}
\eta_n=\max\Big{\{}1-2^{1-n}, \ 1-\frac{2^{2-n}(1-\eta_{n-1})}{2^{2-n}+(n-2)\eta_{n-1}}\Big{\}}.
\end{equation}
If $n=2$, matrix $A$ is diagonalized by only one step of the Jacobi method, which corresponds to one cycle. Therefore, $\eta_2=0$. Hence, from relation~\eqref{3.5} we obtain
$$\eta_4=\max\Big{\{}\frac{7}{8}, \frac{27}{28}\Big{\}}=\frac{27}{28}.$$
\end{proof}

Thus, all orderings from $\mathcal{C}_{sp}^{(n)}$ are convergent.
Let us see how the serial orderings with permutations look like for the case $n=4$.

The relations~\eqref{Ccdef} and~\eqref{Crdef} take the form
$$\mathcal{C}_c^{(4)} = \big{\{} \mathcal{O}\ \big{|} \ \mathcal{O}=(1,2),(\pi_{3}(1),3),(\pi_{3}(2),3),(\pi_{4}(1),4),(\pi_{4}(2),4),(\pi_{4}(3),4) \big{\}},$$
with
$$\pi_3={\small\left(
           \begin{array}{cc}
             1 & 2 \\
             \pi_3(1) & \pi_3(2) \\
           \end{array}
         \right)}, \quad \pi_4={\small\left(
           \begin{array}{ccc}
             1 & 2 & 3 \\
             \pi_4(1) & \pi_4(2) & \pi_4(3) \\
           \end{array}
         \right)},$$
and
$$\mathcal{C}_r^{(4)} = \big{\{}\mathcal{O}\ \big{|} \ \mathcal{O}= (3,4),(2,\tau_{2}(3)),(2,\tau_{2}(4)),(1,\tau_{1}(2)),(1,\tau_{1}(3)),(1,\tau_{1}(4)) \big{\}},$$
with
$$\tau_2={\small\left(
           \begin{array}{cc}
             3 & 4 \\
             \tau_{2}(3) & \tau_{2}(4) \\
           \end{array}
         \right)}, \quad \tau_1={\small\left(
           \begin{array}{ccc}
             2 & 3 & 4 \\
             \tau_{1}(2) & \tau_{1}(3) & \tau_{1}(4) \\
           \end{array}
         \right)},$$
respectively. For the classes of reverse orderings, from~\eqref{Cinvdef} we obtain
\begin{align*}
\overleftarrow{\mathcal{C}_c}^{(4)} & = \big{\{}\mathcal{O} \ \big{|} \ \mathcal{O}= (\pi_{4}(1),4),(\pi_{4}(2),4),(\pi_{4}(3),4), (\pi_{3}(1),3),(\pi_{3}(2),3),(1,2) \big{\}}, \\
\overleftarrow{\mathcal{C}_r}^{(4)} & = \big{\{}\mathcal{O} \ \big{|} \ \mathcal{O}= (1,\tau_{1}(2)),(1,\tau_{1}(3)),(1,\tau_{1}(4)), (2,\tau_{2}(3)),(2,\tau_{2}(4)),(3,4) \big{\}}.
\end{align*}

Let us list the  orderings from $\mathcal{C}_0\cap\mathcal{C}_{sp}^{(4)}$. First we list those from $\mathcal{C}_0\cap\mathcal{C}_c^{(4)}$:
{\scriptsize
\begin{align*}
\mathcal{O}_1 & = (1,2),(1,3),(2,3),(1,4),(2,4),(3,4), \qquad \mathcal{O}_7\,=\, (1,2),(2,3),(1,3),(1,4),(2,4),(3,4), \\
\mathcal{O}_2 & = (1,2),(1,3),(2,3),(1,4),(3,4),(2,4), \qquad \mathcal{O}_8 \,=\, (1,2),(2,3),(1,3),(1,4),(3,4),(2,4), \\
\mathcal{O}_3 & = (1,2),(1,3),(2,3),(2,4),(1,4),(3,4), \qquad \mathcal{O}_9 \,=\, (1,2),(2,3),(1,3),(2,4),(1,4),(3,4), \\
\mathcal{O}_4 & = (1,2),(1,3),(2,3),(2,4),(3,4),(1,4), \qquad \mathcal{O}_{10}= (1,2),(2,3),(1,3),(2,4),(3,4),(1,4), \\
\mathcal{O}_5 & = (1,2),(1,3),(2,3),(3,4),(1,4),(2,4), \qquad \mathcal{O}_{11}= (1,2),(2,3),(1,3),(3,4),(1,4),(2,4), \\
\mathcal{O}_6 & = (1,2),(1,3),(2,3),(3,4),(2,4),(1,4), \qquad \mathcal{O}_{12}= (1,2),(2,3),(1,3),(3,4),(2,4),(1,4).
\end{align*}}
The orderings that belong to $\mathcal{C}_0\cap\overleftarrow{\mathcal{C}_r}^{(4)}$ are the following:
{\scriptsize
\begin{align*}
\mathcal{O}_{13} & = (1,2),(1,3),(1,4),(2,3),(2,4),(3,4), \qquad \mathcal{O}_{15}= (1,2),(1,4),(1,3),(2,3),(2,4),(3,4), \\
\mathcal{O}_{14} & = (1,2),(1,3),(1,4),(2,4),(2,3),(3,4), \qquad \mathcal{O}_{16}= (1,2),(1,4),(1,3),(2,4),(2,3),(3,4).
\end{align*}}
The intersection of the sets $\mathcal{C}_0$ and $\mathcal{C}_r^{(4)}$ is empty because the orderings from $\mathcal{C}_0$ start with pivot pair $(1,2)$ and the orderings from $\mathcal{C}_r^{(4)}$ start with $(3,4)$. Also, $\mathcal{C}_0\cap\overleftarrow{\mathcal{C}_c}^{(4)} = \emptyset$ because the orderings from $\overleftarrow{\mathcal{C}_c}^{(4)}$ start with a pair from the $4$th column.

\subsection{Generalized serial orderings and strategies}\label{sec:generalizedperm}

As has been proved in~\cite{BH_Jac2015}, the set of convergent pivot orderings $\mathcal{C}_{sp}^{(4)}$ can be further generalized to the set $\mathcal{C}_{sg}^{(4)}$ given by
\begin{equation}\label{Csg}
\mathcal{C}_{sg}^{(4)} = \big{\{}\mathcal{O}\in\Ol(\mathbf{P}_4) \ \big{|} \ \mathcal{O}\stackrel{p}{\sim}\mathcal{O'}\stackrel{w}{\sim}\mathcal{O''} \ \ \text{or} \ \ \mathcal{O}\stackrel{w}{\sim}\mathcal{O'}\stackrel{p}{\sim}\mathcal{O''}, \ \ \mathcal{O''}\in\mathcal{C}_{sp}^{(4)}\big{\}}.
\end{equation}
The elements of $\mathcal{C}_{sg}^{(4)}$ are called \emph{generalized serial orderings}.
Proposition~\ref{tm:sg} below, which is a special case of~\cite[Theorem 3.11]{BH_Jac2015}, proves the global convergence of the cyclic Jacobi method defined by orderings from $\mathcal{C}_{sg}^{(4)}$. The constant $\gamma_4$ from relation~\eqref{sg} is the same as $\gamma_4$ in relation~\eqref{sp}.

\begin{Proposition}\label{tm:sg}
Let $A$ be a symmetric $4\times4$ matrix and $\mathcal{O}\in\mathcal{C}_{sg}^{(4)}$. Suppose the chain connecting $\mathcal{O}$ with $\mathcal{O}''\in\mathcal{C}_{sp}^{(4)}$ in~\eqref{Csg} contains $d$ relations of shift-equivalence. Let $A^{[d+1]}$ be obtained from $A$ by applying $d+1$ sweeps of the cyclic Jacobi method defined by the strategy $I_{\mathcal{O}}$.
Then
\begin{equation}\label{sg}
S^2(A^{[d+1]})\leq \eta_4 S^2(A).
\end{equation}
\end{Proposition}

\begin{proof}
The result follows from \cite[Theorem 3.11]{BH_Jac2015} and Proposition~\ref{tm:sp}. In \cite[Theorem 3.11]{BH_Jac2015} one has to use the partition $\pi =(1,1,1,1)$. Proposition~\ref{tm:sp} is used only to calculate $\eta_4$.
\end{proof}

Let us first see how far the set $\mathcal{C}_{sp}^{(4)}$ can be extended without using the permutation equivalence.
In the following relations, for the sake of compactness, we shall mark some orderings from $\Ol(\mathbf{P}_4)$ by $\cdot$.
Moreover, the equivalence relation that transposes $r$-th and $s$-th pair in the ordering will be denoted $\stackrel{r-s}{\sim}$, and the shift-equivalence with the shift of the length $l$, $\stackrel{s(l)}{\sim}$. Finally, in this subsection, the pairs in pivot orderings will be denoted without parenthesis, say $(r,s)$ we will denoted by $\textsf{rs}$. We do this for the conciseness and in order to fit two orderings into one line.

The orderings from $\mathcal{C}_0$ equivalent to some orderings from $\mathcal{C}_c^{(4)}\cup\overleftarrow{\mathcal{C}_r}^{(4)}$ are the following:
{\scriptsize
\begin{align*}
\mathcal{O}_{17} & = \jd,\jt,\jc,\dt,\tc,\dc \stackrel{3-4}{\sim} \ \mathcal{O}_2, \qquad
&\mathcal{O}_{19} & = \jd,\dt,\dc,\jt,\tc,\jc\stackrel{3-4}{\sim} \ \mathcal{O}_{10}, \\
\mathcal{O}_{18} & = \jd,\dt,\dc,\jt,\jc,\tc\stackrel{3-4}{\sim} \ \mathcal{O}_9, \qquad
&\mathcal{O}_{20} & = \jd,\jc,\dc,\jt,\dt,\tc \stackrel{3-4}{\sim} \ \mathcal{O}_{16}.
\end{align*}}
The orderings linked to those from $\mathcal{C}_{sp}^{(4)}$ with one shift-equivalence are as follows:
{\scriptsize
\begin{align*}
\mathcal{O}_{21} & = \jd,\jt,\jc,\tc,\dt,\dc \stackrel{\mathsf{s}(3)}{\sim} \cdot\in\mathcal{C}_r^{(4)}, \qquad
&\mathcal{O}_{36} & = \jd,\dc,\jc,\tc,\jt,\dt \stackrel{\mathsf{s}(1)}{\sim} \cdot\in\overleftarrow{\mathcal{C}_c}^{(4)}, \\
\mathcal{O}_{22} & = \jd,\jt,\jc,\tc,\dc,\dt \stackrel{\mathsf{s}(3)}{\sim} \cdot\in\mathcal{C}_r^{(4)}, \qquad
&\mathcal{O}_{37} & = \jd,\dc,\jc,\tc,\dt,\jt \stackrel{\mathsf{s}(1)}{\sim} \cdot\in\overleftarrow{\mathcal{C}_c}^{(4)}, \\
\mathcal{O}_{23} & = \jd,\jt,\tc,\dt,\dc,\jc \stackrel{\mathsf{s}(2)}{\sim} \cdot\in\mathcal{C}_r^{(4)}, \qquad
&\mathcal{O}_{38} & = \jd,\dc,\tc,\jc,\jt,\dt \stackrel{\mathsf{s}(1)}{\sim} \cdot\in\overleftarrow{\mathcal{C}_c}^{(4)}, \\
\mathcal{O}_{24} & = \jd,\jt,\tc,\dc,\dt,\jc \stackrel{\mathsf{s}(2)}{\sim} \cdot\in\mathcal{C}_r^{(4)}, \qquad
&\mathcal{O}_{39} & = \jd,\dc,\tc,\jc,\dt,\jt \stackrel{\mathsf{s}(1)}{\sim} \cdot\in\overleftarrow{\mathcal{C}_c}^{(4)}, \\
\mathcal{O}_{25} & = \jd,\jc,\jt,\tc,\dt,\dc \stackrel{\mathsf{s}(3)}{\sim} \cdot\in\mathcal{C}_r^{(4)}, \qquad
&\mathcal{O}_{40} & = \jd,\tc,\jc,\dc,\jt,\dt \stackrel{\mathsf{s}(1)}{\sim} \cdot\in\overleftarrow{\mathcal{C}_c}^{(4)}, \\
\mathcal{O}_{26} & = \jd,\jc,\jt,\tc,\dc,\dt \stackrel{\mathsf{s}(3)}{\sim} \cdot\in\mathcal{C}_r^{(4)}, \qquad
&\mathcal{O}_{41} & = \jd,\tc,\jc,\dc,\dt,\jt \stackrel{\mathsf{s}(1)}{\sim} \cdot\in\overleftarrow{\mathcal{C}_c}^{(4)}, \\
\mathcal{O}_{27} & = \jd,\jc,\tc,\dt,\dc,\jt \stackrel{\mathsf{s}(2)}{\sim} \cdot\in\mathcal{C}_r^{(4)}, \qquad
&\mathcal{O}_{42} & = \jd,\tc,\dc,\jc,\jt,\dt \stackrel{\mathsf{s}(1)}{\sim} \cdot\in\overleftarrow{\mathcal{C}_c}^{(4)}, \\
\mathcal{O}_{28} & = \jd,\jc,\tc,\dc,\dt,\jt \stackrel{\mathsf{s}(2)}{\sim} \cdot\in\mathcal{C}_r^{(4)}, \qquad
&\mathcal{O}_{43} & = \jd,\tc,\dc,\jc,\dt,\jt \stackrel{\mathsf{s}(1)}{\sim} \cdot\in\overleftarrow{\mathcal{C}_c}^{(4)}, \\
\mathcal{O}_{29} & = \jd,\tc,\dt,\dc,\jt,\jc \stackrel{\mathsf{s}(1)}{\sim} \cdot\in\mathcal{C}_r^{(4)}, \qquad
&\mathcal{O}_{44} & = \jd,\jt,\dc,\dt,\tc,\jc \stackrel{\mathsf{s}(5)}{\sim} \cdot\in\overleftarrow{\mathcal{C}_r}^{(4)}, \\
\mathcal{O}_{30} & = \jd,\tc,\dt,\dc,\jc,\jt \stackrel{\mathsf{s}(1)}{\sim} \cdot\in\mathcal{C}_r^{(4)}, \qquad
&\mathcal{O}_{45} & = \jd,\jc,\dt,\dc,\tc,\jt \stackrel{\mathsf{s}(5)}{\sim} \cdot\in\overleftarrow{\mathcal{C}_r}^{(4)}, \\
\mathcal{O}_{31} & = \jd,\tc,\dc,\dt,\jt,\jc \stackrel{\mathsf{s}(1)}{\sim} \cdot\in\mathcal{C}_r^{(4)}, \qquad
&\mathcal{O}_{46} & = \jd,\jc,\dc,\dt,\tc,\jt \stackrel{\mathsf{s}(5)}{\sim} \cdot\in\overleftarrow{\mathcal{C}_r}^{(4)}, \\
\mathcal{O}_{32} & = \jd,\tc,\dc,\dt,\jc,\jt \stackrel{\mathsf{s}(1)}{\sim} \cdot\in\mathcal{C}_r^{(4)}, \qquad
&\mathcal{O}_{47} & = \jd,\dt,\dc,\tc,\jt,\jc \stackrel{\mathsf{s}(4)}{\sim} \cdot\in\overleftarrow{\mathcal{C}_r}^{(4)}, \\
\mathcal{O}_{33} & = \jd,\jc,\dc,\tc,\jt,\dt \stackrel{\mathsf{s}(1)}{\sim} \cdot\in\overleftarrow{\mathcal{C}_c}^{(4)}, \qquad
&\mathcal{O}_{48} & = \jd,\dt,\dc,\tc,\jc,\jt \stackrel{\mathsf{s}(4)}{\sim} \cdot\in\overleftarrow{\mathcal{C}_r}^{(4)}, \\
\mathcal{O}_{34} & = \jd,\jc,\dc,\tc,\dt,\jt \stackrel{\mathsf{s}(1)}{\sim} \cdot\in\overleftarrow{\mathcal{C}_c}^{(4)}, \qquad
&\mathcal{O}_{49} & = \jd,\dc,\dt,\tc,\jt,\jc \stackrel{\mathsf{s}(4)}{\sim} \cdot\in\overleftarrow{\mathcal{C}_r}^{(4)}, \\
\mathcal{O}_{35} & = \jd,\jc,\tc,\dc,\jt,\dt \stackrel{\mathsf{s}(1)}{\sim} \cdot\in\overleftarrow{\mathcal{C}_c}^{(4)}, \qquad
&\mathcal{O}_{50} & = \jd,\dc,\dt,\tc,\jc,\jt \stackrel{\mathsf{s}(4)}{\sim} \cdot\in\overleftarrow{\mathcal{C}_r}^{(4)}.
\end{align*}}
Further on, we come to the orderings from $\mathcal{C}_0$ that are weakly equivalent to the ones from $\mathcal{C}_{sp}$. They are listed below:
{\scriptsize
\begin{align*}
\mathcal{O}_{51} & = \jd,\tc,\jt,\dt,\jc,\dc \stackrel{1-2}{\sim} \cdot \stackrel{\mathsf{s}(1)}{\sim} \ \mathcal{O}_1, \qquad
&\mathcal{O}_{58} & = \jd,\dc,\tc,\dt,\jc,\jt \stackrel{4-5}{\sim} \cdot \stackrel{\mathsf{s}(1)}{\sim} \cdot\in\overleftarrow{\mathcal{C}_c}^{(4)}, \\
\mathcal{O}_{52} & = \jd,\tc,\jt,\dt,\dc,\jc \stackrel{1-2}{\sim} \cdot \stackrel{\mathsf{s}(1)}{\sim} \ \mathcal{O}_3, \qquad
&\mathcal{O}_{59} & = \jd,\dt,\jc,\dc,\tc,\jt \stackrel{2-3}{\sim} \cdot \stackrel{\mathsf{s}(5)}{\sim} \cdot\in\overleftarrow{\mathcal{C}_r}^{(4)}, \\
\mathcal{O}_{53} & = \jd,\tc,\dt,\jt,\jc,\dc \stackrel{1-2}{\sim} \cdot \stackrel{\mathsf{s}(1)}{\sim} \ \mathcal{O}_7, \qquad
&\mathcal{O}_{60} & = \jd,\dc,\jt,\dt,\tc,\jc \stackrel{2-3}{\sim} \cdot \stackrel{\mathsf{s}(5)}{\sim} \cdot\in\overleftarrow{\mathcal{C}_r}^{(4)}, \\
\mathcal{O}_{54} & = \jd,\tc,\dt,\jt,\dc,\jc \stackrel{1-2}{\sim} \cdot \stackrel{\mathsf{s}(1)}{\sim} \ \mathcal{O}_9, \qquad
&\mathcal{O}_{61} & = \jd,\tc,\jt,\jc,\dt,\dc \stackrel{1-2}{\sim} \cdot \stackrel{\mathsf{s}(1)}{\sim} \ \mathcal{O}_{13}, \\
\mathcal{O}_{55} & = \jd,\jt,\tc,\dc,\jc,\dt \stackrel{5-6}{\sim} \cdot \stackrel{\mathsf{s}(2)}{\sim} \cdot\in\mathcal{C}_r^{(4)}, \qquad
&\mathcal{O}_{62} & = \jd,\tc,\jt,\jc,\dc,\dt \stackrel{1-2}{\sim} \cdot \stackrel{\mathsf{s}(1)}{\sim} \ \mathcal{O}_{14}, \\
\mathcal{O}_{56} & = \jd,\jc,\tc,\dt,\jt,\dc \stackrel{5-6}{\sim} \cdot \stackrel{\mathsf{s}(2)}{\sim} \cdot\in\mathcal{C}_r^{(4)}, \qquad
&\mathcal{O}_{63} & = \jd,\tc,\jc,\jt,\dt,\dc \stackrel{1-2}{\sim} \cdot \stackrel{\mathsf{s}(1)}{\sim} \ \mathcal{O}_{15}, \\
\mathcal{O}_{57} & = \jd,\jc,\tc,\jt,\dc,\dt \stackrel{4-5}{\sim} \cdot \stackrel{\mathsf{s}(1)}{\sim} \cdot\in\overleftarrow{\mathcal{C}_c}^{(4)}, \qquad
&\mathcal{O}_{64} & = \jd,\tc,\jc,\jt,\dc,\dt \stackrel{1-2}{\sim} \cdot \stackrel{\mathsf{s}(1)}{\sim} \ \mathcal{O}_{16}.
\end{align*}}

Using the relations $\sim$, $\stackrel{s}{\sim}$ and $\stackrel{w}{\sim}$, we have enlarged the set of $16$ orderings belonging to $\mathcal{C}_{sp}^{(4)}$ to the set of $64$ orderings belonging to $\mathcal{C}_{sg}^{(4)}$. To continue, we use the permutation equivalence relation. We will only need four different permutations.
The first one is
$$\mathsf{q}_1=\left(
              \begin{array}{cccc}
                1 & 2 & 3 & 4 \\
                3 & 1 & 2 & 4 \\
              \end{array}
            \right).$$
Applying $\mathsf{p}_1$ we obtain the following convergent orderings:
{\scriptsize
\begin{align*}
\mathcal{O}_{65} & = \jd,\jt,\jc,\dc,\tc,\dt \stackrel{\mathsf{q}_1}{\sim} \cdot \stackrel{\mathsf{s}(5)}{\sim} \ \mathcal{O}_5, \
&\mathcal{O}_{75} & = \jd,\dt,\dc,\jc,\jt,\tc \stackrel{\mathsf{q}_1}{\sim} \cdot \stackrel{\mathsf{s}(3)}{\sim} \cdot\in\mathcal{C}_r^{(4)}, \\
\mathcal{O}_{66} & = \jd,\jt,\dc,\jc,\tc,\dt \stackrel{\mathsf{q}_1}{\sim} \cdot \stackrel{\mathsf{s}(5)}{\sim} \ \mathcal{O}_2, \
&\mathcal{O}_{76} & = \jd,\dc,\jc,\jt,\tc,\dt \stackrel{\mathsf{q}_1}{\sim} \cdot \stackrel{\mathsf{s}(2)}{\sim} \cdot\in\mathcal{C}_r^{(4)}, \\
\mathcal{O}_{67} & = \jd,\jt,\dc,\tc,\jc,\dt \stackrel{\mathsf{q}_1}{\sim} \cdot \stackrel{\mathsf{s}(5)}{\sim} \mathcal{O}_1, \
&\mathcal{O}_{77} & = \jd,\dc,\jc,\dt,\tc,\jt \stackrel{\mathsf{q}_1}{\sim} \cdot \stackrel{3-4}{\sim} \cdot \stackrel{\mathsf{s}(3)}{\sim} \cdot \in\mathcal{C}_r^{(4)}, \\
\mathcal{O}_{68} & = \jd,\jt,\dc,\tc,\dt,\jc \stackrel{\mathsf{q}_1}{\sim} \cdot \stackrel{5-6}{\sim} \cdot \stackrel{\mathsf{s}(5)}{\sim} \ \mathcal{O}_1, \
&\mathcal{O}_{78} & = \jd,\dc,\dt,\jc,\tc,\jt \stackrel{\mathsf{q}_1}{\sim} \cdot \stackrel{\mathsf{s}(3)}{\sim} \cdot\in\mathcal{C}_r^{(4)}, \\
\mathcal{O}_{69} & = \jd,\dc,\jt,\jc,\tc,\dt \stackrel{\mathsf{q}_1}{\sim} \cdot \stackrel{2-3}{\sim} \cdot \stackrel{\mathsf{s}(5)}{\sim} \ \mathcal{O}_2, \
&\mathcal{O}_{79} & = \jd,\jt,\tc,\jc,\dt,\dc \stackrel{\mathsf{q}_1}{\sim} \cdot \stackrel{\mathsf{s}(4)}{\sim} \ \mathcal{O}_{15}, \\
\mathcal{O}_{70} & = \jd,\dc,\jt,\tc,\jc,\dt \stackrel{\mathsf{q}_1}{\sim} \cdot \stackrel{2-3}{\sim} \cdot \stackrel{\mathsf{s}(5)}{\sim} \ \mathcal{O}_1, \
&\mathcal{O}_{80} & = \jd,\jt,\tc,\jc,\dc,\dt \stackrel{\mathsf{q}_1}{\sim} \cdot \stackrel{\mathsf{s}(4)}{\sim} \cdot\in\overleftarrow{\mathcal{C}_r}^{(4)}, \\
\mathcal{O}_{71} & = \jd,\dc,\jt,\tc,\dt,\jc \stackrel{\mathsf{q}_1}{\sim} \cdot \stackrel{2-3}{\sim} \cdot \stackrel{5-6}{\sim} \cdot \stackrel{\mathsf{s}(5)}{\sim} \mathcal{O}_1,
&\mathcal{O}_{81} & = \jd,\dt,\tc,\jt,\jc,\dc \stackrel{\mathsf{q}_1}{\sim} \cdot \stackrel{\mathsf{s}(5)}{\sim} \cdot\in\overleftarrow{\mathcal{C}_r}^{(4)}, \\
\mathcal{O}_{72} & = \jd,\jc,\dt,\jt,\tc,\dc \stackrel{\mathsf{q}_1}{\sim} \cdot \stackrel{2-3}{\sim} \cdot \stackrel{\mathsf{s}(2)}{\sim} \cdot\in\mathcal{C}_r^{(4)}, \
&\mathcal{O}_{82} & = \jd,\dc,\dt,\jt,\tc,\jc \stackrel{\mathsf{q}_1}{\sim} \cdot\in\overleftarrow{\mathcal{C}_r}^{(4)}, \\
\mathcal{O}_{73} & = \jd,\jc,\tc,\jt,\dt,\dc \stackrel{\mathsf{q}_1}{\sim} \cdot \stackrel{\mathsf{s}(1)}{\sim} \cdot\in\mathcal{C}_r^{(4)}, \
&\mathcal{O}_{83} & = \jd,\dc,\tc,\jt,\jc,\dt \stackrel{\mathsf{q}_1}{\sim} \cdot \stackrel{\mathsf{s}(5)}{\sim} \ \mathcal{O}_{14}, \\
\mathcal{O}_{74} & = \jd,\dt,\jc,\jt,\tc,\dc \stackrel{\mathsf{q}_1}{\sim} \cdot \stackrel{\mathsf{s}(2)}{\sim} \cdot\in\mathcal{C}_r^{(4)}, \
&\mathcal{O}_{84} & = \jd,\dc,\tc,\jt,\dt,\jc \stackrel{\mathsf{q}_1}{\sim} \cdot \stackrel{5-6}{\sim} \cdot \stackrel{\mathsf{s}(5)}{\sim} \ \mathcal{O}_{14}.
\end{align*}}
Next, we use the permutation
$$\mathsf{q}_2=\left(
             \begin{array}{cccc}
               1 & 2 & 3 & 4 \\
               1 & 3 & 2 & 4 \\
             \end{array}
           \right).$$
It leads to the following convergent orderings:
{\scriptsize
\begin{align*}
\mathcal{O}_{85} & = \jd,\dt,\jc,\tc,\jt,\dc \stackrel{\mathsf{q}_2}{\sim} \cdot \stackrel{5-6}{\sim} \cdot \stackrel{\mathsf{s}(5)}{\sim} \ \mathcal{O}_1, \
&\mathcal{O}_{90} & = \jd,\jc,\jt,\dt,\tc,\dc \stackrel{\mathsf{q}_2}{\sim} \cdot\in\overleftarrow{\mathcal{C}_r}^{(4)}, \\
\mathcal{O}_{86} & = \jd,\dt,\jc,\tc,\dc,\jt \stackrel{\mathsf{q}_2}{\sim} \cdot \stackrel{\mathsf{s}(5)}{\sim} \ \mathcal{O}_1, \
&\mathcal{O}_{91} & = \jd,\jc,\dt,\tc,\jt,\dc \stackrel{\mathsf{q}_2}{\sim} \cdot \stackrel{5-6}{\sim} \cdot \stackrel{\mathsf{s}(5)}{\sim} \ \mathcal{O}_{13}, \\
\mathcal{O}_{87} & = \jd,\jc,\jt,\dc,\tc,\dt \stackrel{\mathsf{q}_2}{\sim} \cdot \stackrel{\mathsf{s}(3)}{\sim} \cdot\in\mathcal{C}_r^{(4)}, \
&\mathcal{O}_{92} & = \jd,\jc,\dt,\tc,\dc,\jt \stackrel{\mathsf{q}_2}{\sim} \cdot \stackrel{\mathsf{s}(5)}{\sim} \ \mathcal{O}_{13}, \\
\mathcal{O}_{88} & = \jd,\jc,\dc,\jt,\tc,\dt \stackrel{\mathsf{q}_2}{\sim} \cdot \stackrel{3-4}{\sim} \cdot \stackrel{\mathsf{s}(3)}{\sim} \cdot\in\mathcal{C}_r^{(4)},
&\mathcal{O}_{93} & = \jd,\dt,\tc,\jt,\dc,\tc \stackrel{\mathsf{q}_2}{\sim} \cdot \stackrel{4-5}{\sim} \cdot \stackrel{\mathsf{s}(4)}{\sim} \ \mathcal{O}_{15}, \\
\mathcal{O}_{89} & = \jd,\dc,\tc,\dt,\jt,\jc \stackrel{\mathsf{q}_2}{\sim} \cdot \stackrel{\mathsf{s}(1)}{\sim} \cdot\in\mathcal{C}_r^{(4)}, \
&\mathcal{O}_{94} &= \jd,\dt,\tc,\dc,\jt,\jc \stackrel{\mathsf{q}_2}{\sim} \cdot \stackrel{\mathsf{s}(4)}{\sim} \ \mathcal{O}_{15}.
\end{align*}}
We shall also use the permutation equivalence $\stackrel{p}{\sim}$ with
$$\mathsf{q}_3=\left(
                \begin{array}{cccc}
                  1 & 2 & 3 & 4 \\
                  3 & 2 & 1 & 4 \\
                \end{array}
              \right).$$
It yields the following convergent orderings:
{\scriptsize
\begin{align*}
\mathcal{O}_{95} & = \jd,\jc,\dc,\dt,\jt,\tc \stackrel{\mathsf{q}_3}{\sim} \cdot \stackrel{\mathsf{s}(3)}{\sim} \cdot \stackrel{3-4}{\sim} \ \mathcal{O}_2, \
&\mathcal{O}_{100} & = \jd,\dt,\tc,\dc,\jc,\jt \stackrel{\mathsf{q}_3}{\sim} \cdot \stackrel{\mathsf{s}(2)}{\sim} \cdot\in\overleftarrow{\mathcal{C}_c}^{(4)}, \\
\mathcal{O}_{96} & = \jd,\jt,\tc,\dt,\jc,\dc \stackrel{\mathsf{q}_3}{\sim} \cdot \stackrel{\mathsf{s}(4)}{\sim} \cdot\in\mathcal{C}_r^{(4)}, \
&\mathcal{O}_{101} & = \jd,\dc,\jc,\jt,\dt,\tc \stackrel{\mathsf{q}_3}{\sim} \cdot \stackrel{\mathsf{s}(3)}{\sim} \cdot\in\overleftarrow{\mathcal{C}_r}^{(4)}, \\
\mathcal{O}_{97} & = \jd,\dt,\dc,\jc,\tc,\jt \stackrel{\mathsf{q}_3}{\sim} \cdot \stackrel{\mathsf{s}(2)}{\sim} \cdot\in\overleftarrow{\mathcal{C}_c}^{(4)}, \
&\mathcal{O}_{102} & = \jd,\dc,\jc,\dt,\jt,\tc \stackrel{\mathsf{q}_3}{\sim} \cdot \stackrel{\mathsf{s}(3)}{\sim} \ \mathcal{O}_{13}, \\
\mathcal{O}_{98} & = \jd,\dt,\tc,\jc,\jt,\dc \stackrel{\mathsf{q}_3}{\sim} \cdot \stackrel{5-6}{\sim} \cdot \stackrel{\mathsf{s}(4)}{\sim} \cdot\in\overleftarrow{\mathcal{C}_c}^{(4)}, \
&\mathcal{O}_{103} & = \jd,\dc,\dt,\jc,\jt,\tc \stackrel{\mathsf{q}_3}{\sim} \cdot \stackrel{3-4}{\sim} \cdot \stackrel{\mathsf{s}(3)}{\sim} \ \mathcal{O}_{13}. \\
\mathcal{O}_{99} & = \jd,\dt,\tc,\jc,\dc,\jt \stackrel{\mathsf{q}_3}{\sim} \cdot \stackrel{\mathsf{s}(2)}{\sim} \cdot\in\overleftarrow{\mathcal{C}_c}^{(4)},
\end{align*}}
Finally, for one pivot ordering we use the permutation
$\mathsf{q}_4=\left(
                \begin{array}{cccc}
                  1 & 2 & 3 & 4 \\
                  1 & 3 & 4 & 2 \\
                \end{array}
              \right).$
We have
{\scriptsize
$$\mathcal{O}_{104} = \jd,\dc,\dt,\jt,\jc,\tc \stackrel{\mathsf{q}_4}{\sim} \cdot \stackrel{\mathsf{s}(3)}{\sim} \cdot \stackrel{3-4}{\sim} \cdot\in\overleftarrow{\mathcal{C}_r}^{(4)}.$$}

There are $16$ orderings left in $\mathcal{C}_0\big{\backslash}\mathcal{C}_{sg}^{(4)}$ that cannot be linked to any of the orderings from $\mathcal{C}_{sg}^{(4)}$ by using the equivalence relations $\sim$, $\stackrel{s}{\sim}$, $\stackrel{w}{\sim}$ or $\stackrel{p}{\sim}$.
These are the parallel orderings.

\subsection{Parallel strategies}\label{sec:parallel}

Behavior of the Jacobi method under parallel pivot strategies has been studied in~\cite{BH_4x4pa2015}. Parallel pivot strategies are interesting especially because, unlike the serial ones, they can force the Jacobi method to be very fast or very slow within one cycle.
As a representative of parallel orderings for $n=4$ we take the pivot ordering
\begin{equation}\label{Opar}
\mathcal{O}_{\text{par}}=(1,3),(2,4),(1,4),(2,3),(1,2),(3,4).
\end{equation}
Its ``parallelism'' manifests in the fact that pairs
\begin{equation}\label{Opar_parovi}
\{(1,3),(2,4)\}, \ \{(1,4),(2,3)\} \ \text{and} \ \{(1,2),(3,4)\}
\end{equation}
commute. Therefore, the corresponding Jacobi transformations commute and can be computed and applied independently.
Instead of six sequential steps needed to complete one cycle, only three parallel steps are used.

If $\mathcal{O}\sim\mathcal{O}_{\text{par}}$, it is easy to see that the only allowed transpositions are those of the pairs within the curly braces from~\eqref{Opar_parovi}.
Thus, after each parallel step, the Jacobi method defined by such $I_{\mathcal{O}}$ generates the same matrix as the Jacobi method defined by $I_{\mathcal{O}_{\text{par}}}$.

Another representative of parallel orderings is
\begin{equation}\label{Opar2}
\mathcal{O}_{\text{par}}'= (1,4),(2,3),(1,3),(2,4),(1,2),(3,4).
\end{equation}
It is easily seen that $\mathcal{O}_{\text{par}}'$ is not weakly equivalent to $\mathcal{O}_{\text{par}}$. However, these two orderings are permutation equivalent, in particular, we have
$$\mathcal{O}_{\text{par}}'= \mathcal{O}_{\text{par}}(\mathsf{p}_{\text{par}}),\qquad
\mathsf{p}_{\text{par}}=\left(\begin{array}{cccc}
                  1 & 2 & 3 & 4 \\
                  1 & 2 & 4 & 3 \\
                \end{array}\right).$$

If $\mathcal{O}\stackrel{p}{\sim}\mathcal{O}_{\text{par}}$, then $\mathcal{O}$ results from $\mathcal{O}_{\text{par}}$ or $\mathcal{O}_{\text{par}}'$ by applying an appropriate $\sim$ and afterwards an appropriate $\stackrel{s}{\sim}$ with shift of length $0$, $2$ or $4$, or vice versa (first $\stackrel{s}{\sim}$ then $\sim$). The same is true for $\mathcal{O}\stackrel{p}{\sim}\mathcal{O}_{\text{par}}'$. This follows from the fact that permutations are bijections and commuting pairs are mapped into commuting pairs.

For the Jacobi method defined by ordering $\mathcal{O}_{\text{par}}$ from~\eqref{Opar} or $\mathcal{O}_{\text{par}}'$ from~\eqref{Opar2}, Proposition~\ref{tm:par} holds. Note that the condition $a_{12}=0$ and $a_{34}=0$ is not a restriction when we are considering the global convergence. If that condition does not hold for the original matrix $A$ it will hold at the beginning of each cycle of the Jacobi method, starting with the second one.

\begin{Proposition}\label{tm:par}
Let $A$ be a symmetric $4\times4$ matrix such that $a_{12}=0$, $a_{34}=0$, and let $A^{[2]}$ be obtained by applying two sweeps of the Jacobi method under the strategy $I_{\mathcal{O}_{\text{par}}}$ ($I_{\mathcal{O}_{\text{par}}'}$) to $A$. Then
\begin{equation}\label{eqpar}
S(A^{[2]})\leq \gamma_4'S(A), \quad \gamma_4' = 1-10^{-5}.
\end{equation}
Estimate~\eqref{eqpar} holds if $\mathcal{O}_{\text{par}}$ ($\mathcal{O}_{\text{par}}'$) is replaced by any $\mathcal{O}$ such that $\mathcal{O}\sim\mathcal{O}_{\text{par}}$ ($\mathcal{O}\sim\mathcal{O}_{\text{par}}'$).
\end{Proposition}

\begin{proof}
The first assertion for $I_{\mathcal{O}_{\text{par}}}$ is~\cite[Theorem 3.4]{BH_4x4pa2015}. The second one for $\mathcal{O}\sim \mathcal{O}_{\text{par}}$ follows from the fact that the Jacobi methods defined by the strategies $I_{\mathcal{O}_{\text{par}}}$ and $I_{\mathcal{O}}\sim I_{\mathcal{O}_{\text{par}}}$ generate the same matrix $A^{[2]}$.

To prove the assertion for $I_{\mathcal{O}_{\text{par}}'}$, let $\mathsf{P}=[e_1,e_2,e_4,e_3]$, where $[e_1,e_2,e_3,e_4]=I$. Let $(A^{(k)},k\geq0)$ and $(\tilde{A}^{(k)},k\geq0)$ be two sequences of matrices generated by the Jacobi method defined by the strategies $I_{\mathcal{O_{\text{par}}}}$ and $I_{\mathcal{O_{\text{par}}}'}$, respectively.
From~\cite[Proposition 3.1]{BH_4x4pa2015} we have
$$\tilde{A}^{(2r)}=\mathsf{P}^TA^{(2r)}\mathsf{P}, \quad r\geq0.$$
Therefore, $S(\tilde{A}^{[0]})=S(A^{[0]})$ and $S(\tilde{A}^{[2]})=S(A^{[2]})$, which proves the relation~\eqref{eqpar} for $I_{\mathcal{O}_{\text{par}}'}$.
The final assertion for $\mathcal{O}\sim \mathcal{O}_{\text{par}}'$ is obtained by the same argument as for $\mathcal{O}\sim\mathcal{O}_{\text{par}}$, but $I_{\mathcal{O}_{\text{par}}'}$ is used instead of $I_{\mathcal{O}_{\text{par}}}$.
\end{proof}

The following corollary uses the same $\gamma_4'$.

\begin{Corollary}\label{cor:par}
Let $\mathcal{O},\mathcal{O}',\mathcal{O}''\in\Ol(\mathbf{P}_4)$ be such that either $\mathcal{O}\stackrel{s}{\sim}\mathcal{O}'\sim\mathcal{O}_{\text{par}}$ or  $\mathcal{O}\sim\mathcal{O}''\stackrel{s}{\sim}\mathcal{O}_{\text{par}}$. Let $A$ be a symmetric $4\times4$ matrix and let $A^{[t]}$ be obtained after applying $t$, $t\geq1$, sweeps of the Jacobi method under the strategy $I_{\mathcal{O}}$ to $A$. Then
\begin{equation}\label{eqcor}
S(A^{[t+3]})\leq\gamma_4'S(A^{[t]})).
\end{equation}
Moreover, if $t=0$, then the relation~\eqref{eqcor} holds if the length of the shift is not one or if $a_{12}=a_{34}=0$.

The same assertion holds provided that $\mathcal{O}_{\text{par}}$ is replaced by $\mathcal{O}_{\text{par}}'$.
\end{Corollary}

\begin{proof}
First, consider the case $\mathcal{O}\stackrel{s}{\sim}\mathcal{O}_{\text{par}}$.
The sequence of pivot pairs $[\mathcal{O},\mathcal{O},\mathcal{O}]$ contains the sequence $[\mathcal{O}_{\text{par}},\mathcal{O}_{\text{par}}]$. Since $t\geq1$, Proposition~\ref{tm:par} can be applied to $A^{[t]}$ and assertion~\eqref{eqcor} follows.

If the length of shift is not one, then the sequence of the pairs $(1,2),(3,4)$ precedes the sequence $[\mathcal{O}_{\text{par}},\mathcal{O}_{\text{par}}]$. Elements at the positions $(1,2)$ and $(3,4)$ will both be zero prior to the application of Proposition~\ref{tm:par}. Thus, relation~\eqref{eqcor} holds for $t=0$ as well.
If the length of shift is one, only the pair $(3,4)$ will precede the sequence $[\mathcal{O}_{\text{par}},\mathcal{O}_{\text{par}}]$. Therefore, in order to apply Proposition~\ref{tm:par}, we should require $a_{12}=0$.

If $\mathcal{O}\stackrel{s}{\sim}\mathcal{O}'\sim\mathcal{O}_{\text{par}}$,
the reasoning is the same, where the role of $\mathcal{O}_{\text{par}}$ is now played by $\mathcal{O}'$. Sequence $\mathcal{O}'$ can differ from $\mathcal{O}_{\text{par}}$ in at most three transpositions of the pairs from~\eqref{Opar_parovi}. Therefore, if the length of the shift is one, one of the pairs $(1,2),(3,4)$ will precede the sequence $[\mathcal{O}_{\text{par}},\mathcal{O}_{\text{par}}]$ and we should require $a_{12}=a_{34}=0$.

If $\mathcal{O}\sim\mathcal{O}''\stackrel{s}{\sim}\mathcal{O}_{\text{par}}$,
the reasoning is similar, only the role of $\mathcal{O}$ is now played by $\mathcal{O}''$. Since $\mathcal{O}\sim\mathcal{O}''$ the sequence $\mathcal{O}$ can differ from $\mathcal{O}''$ in at most three transpositions of the pairs from~\eqref{Opar_parovi}. If the length of the shift is one, only one of the pairs $(1,2),(3,4)$ will precede the sequence $[\mathcal{O}_{\text{par}},\mathcal{O}_{\text{par}}]$ and we should require $a_{12}=a_{34}=0$.

The assertions which use $\mathcal{O}_{\text{par}}'$ are proved in the same way. Every appearance of $\mathcal{O}_{\text{par}}$ is replaced by $\mathcal{O}_{\text{par}}'$.
\end{proof}

There are $16$ orderings from $\mathcal{C}_0$ that can be fully parallelized. These orderings are exactly those from $\mathcal{C}_0\big{\backslash}\mathcal{C}_{sg}^{(4)}$.
Using the equivalence and shift-equivalence relations (i.e.\ the weak equivalence), we can link the following eight orderings from $\mathcal{C}_0$ with $\mathcal{O}_{\text{par}}$:
{\scriptsize
\begin{align*}
\mathcal{O}_{105} & = (1,2),(3,4),(1,3),(2,4),(1,4),(2,3) \stackrel{\mathsf{s}(2)}{\sim} \ \mathcal{O}_{\text{par}}, \\
\mathcal{O}_{106} & = (1,2),(1,3),(2,4),(1,4),(2,3),(3,4) \stackrel{\mathsf{s}(1)}{\sim} \cdot \stackrel{5-6}{\sim} \ \mathcal{O}_{\text{par}}, \\
\mathcal{O}_{107} & = (1,2),(1,3),(2,4),(2,3),(1,4),(3,4) \stackrel{\mathsf{s}(1)}{\sim} \cdot \stackrel{3-4}{\sim} \cdot \stackrel{5-6}{\sim} \ \mathcal{O}_{\text{par}}, \\
\mathcal{O}_{108} & = (1,2),(2,4),(1,3),(1,4),(2,3),(3,4) \stackrel{\mathsf{s}(1)}{\sim} \cdot \stackrel{1-2}{\sim} \cdot \stackrel{5-6}{\sim} \ \mathcal{O}_{\text{par}}, \\
\mathcal{O}_{109} & = (1,2),(2,4),(1,3),(2,3),(1,4),(3,4) \stackrel{\mathsf{s}(1)}{\sim} \cdot \stackrel{1-2}{\sim} \cdot \stackrel{3-4}{\sim} \cdot \stackrel{5-6}{\sim} \ \mathcal{O}_{\text{par}}, \\
\mathcal{O}_{110} & = (1,2),(3,4),(1,3),(2,4),(2,3),(1,4) \stackrel{\mathsf{s}(2)}{\sim} \cdot \stackrel{3-4}{\sim} \ \mathcal{O}_{\text{par}}, \\
\mathcal{O}_{111} & = (1,2),(3,4),(2,4),(1,3),(1,4),(2,3) \stackrel{\mathsf{s}(2)}{\sim} \cdot \stackrel{1-2}{\sim} \ \mathcal{O}_{\text{par}}, \\
\mathcal{O}_{112} & = (1,2),(3,4),(2,4),(1,3),(2,3),(1,4) \stackrel{\mathsf{s}(2)}{\sim} \cdot \stackrel{1-2}{\sim} \cdot \stackrel{3-4}{\sim} \ \mathcal{O}_{\text{par}}.
\end{align*}}
In the same manner the other eight orderings from $\mathcal{C}_0$ can be linked with $\mathcal{O}_{\text{par}}'$:
{\scriptsize
\begin{align*}
\mathcal{O}_{113} & = (1,2),(1,4),(2,3),(1,3),(2,4),(3,4) \stackrel{\mathsf{s}(1)}{\sim} \cdot \stackrel{5-6}{\sim} \ \mathcal{O}_{\text{par}}', \\
\mathcal{O}_{114} & = (1,2),(1,4),(2,3),(2,4),(1,3),(3,4) \stackrel{\mathsf{s}(1)}{\sim} \cdot \stackrel{3-4}{\sim} \cdot \stackrel{5-6}{\sim} \ \mathcal{O}_{\text{par}}', \\
\mathcal{O}_{115} & = (1,2),(2,3),(1,4),(1,3),(2,4),(3,4) \stackrel{\mathsf{s}(1)}{\sim} \cdot \stackrel{1-2}{\sim} \cdot \stackrel{5-6}{\sim} \ \mathcal{O}_{\text{par}}', \\
\mathcal{O}_{116} & = (1,2),(2,3),(1,4),(2,4),(1,3),(3,4) \stackrel{\mathsf{s}(1)}{\sim} \cdot \stackrel{1-2}{\sim} \cdot \stackrel{3-4}{\sim} \cdot \stackrel{5-6}{\sim} \ \mathcal{O}_{\text{par}}', \\
\mathcal{O}_{117} & = (1,2),(3,4),(1,4),(2,3),(1,3),(2,4) \stackrel{\mathsf{s}(2)}{\sim} \ \mathcal{O}_{\text{par}}', \\
\mathcal{O}_{118} & = (1,2),(3,4),(1,4),(2,3),(2,4),(1,3) \stackrel{\mathsf{s}(2)}{\sim} \cdot \stackrel{3-4}{\sim} \ \mathcal{O}_{\text{par}}', \\
\mathcal{O}_{119} & = (1,2),(3,4),(2,3),(1,4),(1,3),(2,4) \stackrel{\mathsf{s}(2)}{\sim} \cdot \stackrel{1-2}{\sim} \ \mathcal{O}_{\text{par}}', \\
\mathcal{O}_{120} & = (1,2),(3,4),(2,3),(1,4),(2,4),(1,3) \stackrel{\mathsf{s}(2)}{\sim} \cdot \stackrel{1-2}{\sim} \cdot \stackrel{3-4}{\sim} \ \mathcal{O}_{\text{par}}'.
\end{align*}}

This shows that all $120$ orderings from $\mathcal{C}_0$ are convergent.

\section{Proof of Theorem~\ref{tm:main}}\label{sec:proof}

To prove Theorem~\ref{tm:main} we use the results from sections~\ref{sec:strategies} and
\ref{sec:strategies3}.
We know that there are $720$ cyclic orderings in $\Ol(\mathbf{P}_4)$ and $120$ of them are in $\mathcal{C}_0$.
Each ordering from $\Ol(\mathbf{P}_4)$ either belongs to $\mathcal{C}_{sg}^{(4)}$ or is weakly equivalent to either $\mathcal{O}_{par}$ or $\mathcal{O}_{par}'$.

Let $A$ be an arbitrary symmetric $4\times4$ matrix. By $A^{[t]}$ we denote the matrix obtained from $A$ after $t$ sweeps of the Jacobi method.
Recall that $\gamma_4'=1-10^{-5}$ and $\eta_4=\frac{27}{28}$. Set $\gamma_4=\sqrt{\eta_4}=\frac{3\sqrt{21}}{14}<1$.

First we deal with the pivot orderings $\mathcal{O}_i$, $1\leq i\leq120$, from $\mathcal{C}_0$ that have been analyzed in Section~\ref{sec:strategies3}.
For $i\in\{1,\ldots,16\}$, we have $\mathcal{O}_i\in\mathcal{C}_{sp}^{(4)}$ and Proposition~\ref{tm:sp} implies
\begin{equation}\label{proof:sp}
S(A^{[1]})\leq\gamma_4 S(A).
\end{equation}
For $i\in\{17,\ldots,104\}$, we use Proposition~\ref{tm:sg}.
If $i\in\{17,18,19,20,82,90\}$, then there are no shift-equivalences linking $\mathcal{O}_i$ to an ordering from $\mathcal{C}_{sp}^{(4)}$. Therefore, for the corresponding pivot strategies relation~\eqref{proof:sp} holds.
Otherwise, there is one shift-equivalence in the corresponding links and we have
\begin{equation}\label{proof:sg}
S(A^{[2]})\leq\gamma_4 S(A).
\end{equation}

For $i\in\{105,\ldots,120\}$, the orderings $\mathcal{O}_i$ are linked to the ordering $\mathcal{O}_{par}$ from~\eqref{Opar} or to $\mathcal{O}_{par}'$ from~\eqref{Opar2}
with one shift-equivalence.
For $i\in\{105,110,$ $111,112\}$ the length of the shift is two and from Proposition~\ref{tm:par} it follows
\begin{equation}\label{proof:par}
S(A^{[3]})\leq\gamma_4' S(A).
\end{equation}
For $i\in\{117,\ldots,120\}$ the length of the shift is also two and the relation\eqref{proof:par} follows from the same proposition.
If $i\in\{106,\ldots,109,113,\ldots,116\}$ the length of the shift is one, so we use Corollary~\ref{cor:par} to obtain
\begin{equation}\label{eqtm}
S(A^{[t+3]})\leq\gamma_4' S(A^{[t]}), \quad t\geq1.
\end{equation}

Assume now that $\mathcal{O}\in\Ol(\mathbf{P}_4)\big{\backslash}\mathcal{C}_0$. Then $\mathcal{O}$ is shift-equivalent to some $\mathcal{O}_i\in\mathcal{C}_0$.
If $\mathcal{O}$ is shift-equivalent to $\mathcal{O}_i$, $1\leq i\leq104$, then $\mathcal{O}$ can be linked to some $\mathcal{O}''\in\mathcal{C}_{sp}$ with at most two shift-equivalences. Hence, applying Proposition~\ref{tm:sg} we obtain
$$S(A^{[3]})\leq\gamma_4 S(A).$$

Finally, let $\mathcal{O}\in\Ol(\mathbf{P}_4)\big{\backslash}\mathcal{C}_0$ be shift-equivalent to some $\mathcal{O}_i\in\mathcal{C}_0$, $105\leq i\leq120$. Then, one of the two possibilities holds.
\begin{itemize}
\item[(i)] $\mathcal{O}\stackrel{s}{\sim}\mathcal{O}_i$, $i\in\{105,\ldots,112\}$, and there is $\mathcal{O}'\in\Ol(\mathbf{P}_4)$ such that     $$\mathcal{O}\stackrel{s}{\sim}\mathcal{O}_i\stackrel{s}{\sim}\mathcal{O}'\sim\mathcal{O}_{par}.$$
    By the transitivity property of $\stackrel{s}{\sim}$ we have $\mathcal{O}\stackrel{s}{\sim}\mathcal{O}'\sim\mathcal{O}_{par}$.
\item[(ii)] $\mathcal{O}\stackrel{s}{\sim}\mathcal{O}_i$, $i\in\{113,\ldots,120\}$, and there is $\mathcal{O}''\in\Ol(\mathbf{P}_4)$ such that $$\mathcal{O}\stackrel{s}{\sim}\mathcal{O}_i\stackrel{s}{\sim}\mathcal{O}'\sim\mathcal{O}_{par}, \quad \text{thus} \quad \mathcal{O}\stackrel{s}{\sim}\mathcal{O}''\sim\mathcal{O}_{par}'.$$
\end{itemize}
In both cases we use Corollary~\ref{cor:par} to obtain relation~\eqref{eqtm}.

Summing up, for $n=4$ and for every cyclic Jacobi method it holds
$$S^2(A^{[t+3]})\leq\gamma S^2(A^{[t]}), \quad t\geq1, \quad \gamma=1-10^{-5}.$$
For certain cyclic Jacobi methods better estimates hold. They are given by the relations~\eqref{proof:sp} --~\eqref{proof:par}.

\begin{Remark}
The results of this paper can be generalized to the block methods where an $n\times n$ matrix is partitioned in $\frac{n}{4}\times\frac{n}{4}$ blocks. For all cyclic strategies apart from those that are weakly equivalent to the parallel strategies from Subsection~\ref{sec:parallel}, one should only check that the assumptions from~\cite{BH_Jac2015} are met (see~\cite{BH_Jac2015}, Theorem 3.11). For the strategies weakly equivalent to the parallel ones, Proposition~\ref{tm:par} and Corollary~\ref{cor:par} should be established for the block Jacobi method, with the constant $\gamma_n'<1$, following the ideas from~\cite{BH_4x4pa2015}.
\end{Remark}

\section{An application to $J$-Jacobi method}\label{sec:appl}

The results obtained in this paper have applications to different Jacobi-type methods. Here, we apply it to the $J$-Jacobi method from~\cite{Ves-JJac93,HSS14}.

The $J$-Jacobi method is an iterative method for solving the generalized eigenvalue problem $Ax=\lambda Jx$, $x\neq0$, where $A$ is symmetric, $J$ is diagonal matrix of signs, and the pair $(A,J)$ is definite. Its main application is as a part of the highly accurate composite algorithm for the eigenvalue problem $Hx=\lambda x$, $x\ne0$, with indefinite symmetric matrix $H$. Roughly speaking, after decomposing $H$ into $LJL^T$ by the Bunch-Parlett algorithm, the eigenvalue problem for $H$ is transformed into $Ax=\lambda Jx$, where $A=L^TL$ and $J$ is the diagonal matrix of signs which reveals the inertia of $H$ (see~\cite{HSS14} for the details). The $J$-Jacobi method, which has been proposed by Veseli\'{c}~\cite{Ves-JJac93}, applies $J$-orthogonal transformation matrices to $A$ in order to diagonalize it. Once $A$ is diagonalized by some $J$-orthogonal matrix $F$, that is $F^TAF=\Lambda_A$, $F^TJF=J$, the eigenvalues (eigenvectors) of $H$ are obtained as the diagonal elements of $J\Lambda_A$ (the columns of $L^{-1}JF\Lambda_A^{1/2}$). In practical computation, one-sided version of the $J$-Jacobi method is applied to $L$, but here we consider the two-sided method since the global convergence of the one-sided method reduces to the global convergence of the two-sided method.

Each iteration has the form
\begin{equation}\label{JJacobi}
A^{(k+1)}=F_k^TA^{(k)}F_k, \quad k\geq0, \qquad A^{(0)}=A,
\end{equation}
where $A$ is symmetric positive definite.
The transformation matrices are $J$-orthogonal, which means that $F_k^TJF_k=J$, $k\geq 0$. In particular, each $F_k$ is either a ``trigonometric Jacobi rotation'' or a ``hyperbolic Jacobi rotation'' in which hyperbolic sine and cosine appear. The corresponding step will be called trigonometric or hyperbolic.

The $J$-Jacobi method is said to be globally convergent if, for any initial $A$, the sequence of matrices $(A^{(k)},k\geq0)$ generated by the iterative process~\eqref{JJacobi} converges to some diagonal matrix. Again, it is sufficient to prove that $S(A^{(k)})\rightarrow 0$ as $k\rightarrow\infty$. Then the diagonal elements will converge (see the proof of \cite[Theorem~3.7]{HSS10}). We note that the quadratic convergence for the serial method has been proved in \cite{drm+har-93}.

We are concerned with $J$-symmetric $4\times4$ matrices. Such matrices appear naturally as pivot submatrices in the block $J$-Jacobi method with $2\times2$ blocks (or all column widths equal to $2$).
Hence, we shall assume that $J=\text{diag}(1,1,-1,-1)$. (The cases $J=\text{diag}(1,-1,-1,-1)$ and $J=\text{diag}(1,1,1,-1)$ require somewhat different approach.)

\begin{Corollary}\label{tm:JJac}
Let $A$ be a symmetric $4\times4$ matrix and $J=\text{diag}(1,1,-1,-1)$. Then every cyclic $J$-Jacobi method for the pair $(A,J)$ is globally convergent.
\end{Corollary}

\begin{proof}
Denote the pivot strategy by $I_{\mathcal{O}}$. We distinguish two cases: (i) $\mathcal{O}$ is generalized serial and (ii) $\mathcal{O}$ is weakly equivalent to $\mathcal{O}_{\text{par}}$ or $\mathcal{O}_{\text{par}'}$.

For the case $(i)$ the proof follows directly from~\cite[Theorem~5.4]{BH_Jac2015}, by taking $n=4$, $\nu=2$ and $\pi =(1,1,1,1)$.
Proposition~\ref{tm:sg} (or~\cite[Theorem 3.11]{BH_Jac2015}) is used to ensure that every generalized serial ordering is convergent.

For the case $(ii)$ the proof is different. It uses two facts.

The first one is that the sequence $[\mathcal{O},\mathcal{O}]$ contains the sequence $\tilde{\mathcal{O}}$ which is equivalent to either $(1,2),(3,4),(1,3),(2,4),(1,4),(2,3)$ or $(1,2),(3,4),(1,4),$ $(2,3),(1,3),(2,4)$. We shall consider the first case for $\tilde{\mathcal{O}}$; the proof for the second case is quite similar.

The second fact reads: the ``hyperbolic angle'' $\theta_k$ of the hyperbolic rotation $F_k$ tends to $0$ as $k$ increases \cite[Lemma~2.2]{Ves-JJac93}. That result holds for general $n$ and for any pivot strategy.
This implies that for every $\varepsilon>0$ there is $k_0\in\mathbb{N}$ such that, for $k\geq k_0$,
the moduli of all pivot elements $|a_{i(k)j(k)}^{(k)}|$ associated with hyperbolic steps and all $|\tanh\theta_k|$ are smaller than $\varepsilon/2$.
Since the pivot strategy $I_{\tilde{\mathcal{O}}}$ is parallel, we can write this condition in the form: for every $\varepsilon>0$ there is $r_0\in \mathbb{N}$ such that,
for every $r\geq r_0$,
\begin{align}
\label{piv_par} \max \left\{ (a_{13}^{(6r+2)})^2+(a_{24}^{(6r+2)})^2, \ (a_{14}^{(6r+4)})^2+(a_{23}^{(6r+4)})^2 \right\} & < \frac{\varepsilon^2}{2}, \\
\label{hyp_angl} \max_{1\leq l\leq 2} \left\{ \tanh^2\theta_{6r+2l}+\tanh^2\theta_{6r+2l+1} \right\} & < \frac{\varepsilon^2}{2}.
\end{align}

We shall prove that the relations (\ref{piv_par}) and (\ref{hyp_angl}) imply $\lim_{r\rightarrow\infty}S(A^{(2r)})=0$.
It is enough to prove that for every $\varepsilon$,
\begin{equation}\label{eps01}
0<\varepsilon < 0.1,
\end{equation}
there is $r_0$ such that $S(A^{(2r)})<\varepsilon$, for all $r\geq r_0$.

Let us simplify the notation. For given $r$, $r\geq r_0$, denote
\begin{align*}
& A^{(6r+2)}=(a_{rs}), \ \ A^{(6r+4)}=(a_{rs}'), \ \ A^{(6r+6)}=(a_{rs}''), \ \ A^{(6r+8)}=(a_{rs}'''), \\
& \ch_l=\cosh\theta_{6r+l+1}, \ \sh_l=\sinh\theta_{6r+l+1}, \ \th_l=\tanh\theta_{6r+l+1}, \quad 1\leq l\leq 4.
\end{align*}
Then, we write
{\small
\begin{align*}
\left[\begin{array}{cccc}
a_{11}' & a_{12}' & 0 & a_{14}' \\
a_{12}' & a_{22}' & a_{23}' & 0 \\
0 & a_{23}' & a_{33}' & a_{34}' \\
a_{14}' & 0 & a_{34}' & a_{44}'
\end{array}\right]
& =
\left[\begin{array}{cccc}
\ch_1 & 0 & \sh_1 & 0 \\
0 & \ch_2 & 0 & \sh_2 \\
\sh_1 & 0 & \ch_1 & 0 \\
0 & \sh_2 & 0 & \ch_2
\end{array}\right]
\left[\begin{array}{cccc}
a_{11} & 0 & a_{13} & a_{14} \\
0 & a_{22} & a_{23} & a_{24} \\
a_{14} & a_{24} & a_{44} & 0 \\
a_{13} & a_{23} & 0 & a_{33}
\end{array}\right]
\left[\begin{array}{cccc}
\ch_1 & 0 & \sh_1 & 0 \\
0 & \ch_2 & 0 & \sh_2 \\
\sh_1 & 0 & \ch_1 & 0 \\
0 & \sh_2 & 0 & \ch_2
\end{array}\right], \\
\left[\begin{array}{cccc}
a_{11}''&a_{12}''&a_{13}''&0\\
a_{12}''&a_{22}''&0&a_{24}''\\
a_{33}''&0&a_{33}''&a_{34}''\\
0&a_{24}''&a_{34}''&a_{44}''
\end{array}\right]
& =
\left[\begin{array}{cccc}
ch_3 & 0 & 0 & \sh_3 \\
0 & \ch_4 & \sh_4 & 0 \\
0 & \sh_4 & \ch_4 & 0 \\
\sh_3 & 0 & 0 & \ch_3
\end{array}\right]
\left[\begin{array}{cccc}
a_{11}' & a_{12}' & 0 & a_{14}' \\
a_{12}' & a_{22}' & a_{23}' & 0 \\
0 & a_{23}' & a_{33}' & a_{34}' \\
a_{14}' & 0 & a_{34}' & a_{44}'
\end{array}\right]
\left[\begin{array}{cccc}
ch_3 & 0 & 0 & \sh_3 \\
0 & \ch_4 & \sh_4 & 0 \\
0 & \sh_4 & \ch_4 & 0 \\
\sh_3 & 0 & 0 & \ch_3
\end{array}\right].
\end{align*}}
It holds
\begin{align*}
a_{12}' & = \sh_1 \ch_2 a_{23} + \ch_1 \sh_2 a_{14}, \qquad
a_{14}' = \sh_1 \sh_2 a_{23} + \ch_1 \ch_2 a_{14}, \\
a_{34}' & = \ch_1 \sh_2 a_{23} + \sh_1 \ch_2 a_{14}, \qquad
a_{23}' = \ch_1 \ch_2 a_{23} + \sh_1 \sh_2 a_{14}, \\
a_{12}'' & = \sh_3 \ch_4 a_{12}' + \sh_3 \sh_4 a_{34}', \qquad
a_{13}'' = \ch_3 \sh_4 a_{12}' + \sh_3 \ch_4 a_{34}', \\
a_{34}'' & = \sh_3 \sh_4 a_{12}' + \ch_3 \ch_4 a_{34}', \qquad
a_{24}'' = \sh_3 \ch_4 a_{12}' + \ch_3 \sh_4 a_{34}'.
\end{align*}

It is easy to to prove that
\begin{align*}
1\leq (\ch_1\ch_2-|\sh_1\sh_2|)^2 \leq \frac{(a_{14}')^2+(a_{23}')^2}{a_{14}^2+a_{23}^2} & \leq (\ch_1\ch_2+|\sh_1\sh_2|)^2, \\
\frac{(a_{12}')^2+(a_{34}')^2}{a_{14}^2+a_{23}^2} & \leq (\ch_1 |\sh_2|+ \ch_2 |\sh_1|)^2.
\end{align*}
Note that $a_{14}^2+a_{23}^2=0$ would imply $S(A^{(6r+4)})=0$.
Using the relations~\eqref{hyp_angl} and~\eqref{eps01} we obtain
\begin{align*}
(\ch_1\ch_2+|\sh_1\sh_2|)^2 & = \frac{(1+|\th_1\th_2|)^2}{(1-\th^2_1)(1-\th^2_2)} < \frac{(1+\varepsilon^2/4)^2}{(1-\varepsilon^2/4)^2} < 1+1.0052\varepsilon^2, \\
(\ch_1 |\sh_2|+ \ch_2 |\sh_1|)^2 & = \frac{(|\th_1|+ |\th_2|)^2}{(1-\th^2_1)(1-\th^2_2)} < \frac{\varepsilon^2}{1-\varepsilon^2/4} < 1.00251\varepsilon^2.
\end{align*}
Hence,
\begin{align}
\label{aa1} a_{14}^2+a_{23}^2 \leq (a_{14}')^2 + (a_{23}')^2 & < (1+1.0052\varepsilon^2)(a_{14}^2+a_{23}^2), \\
\label{aa2} (a_{12}')^2 + (a_{34}')^2 & < 1.00251\varepsilon^2(a_{14}^2+a_{23}^2).
\end{align}
Moreover, note that $(a_{12}')^2+(a_{34}')^2=0$ would imply $S(A^{(6r+6)})=0$ and in the same way we have
\begin{align}
\label{aa3} (a_{12}')^2+(a_{34}')^2 \leq (a_{12}'')^2+(a_{34}'')^2 & < (1+1.0052\varepsilon^2 ) \big{(}(a_{12}')^2+(a_{34}')^2\big{)}, \\
\label{aa4} (a_{13}'')^2+(a_{24}'')^2 & < 1.00251\varepsilon^2 \big{(}(a_{12}')^2+(a_{34}')^2\big{)}.
\end{align}

We use the relations~\eqref{piv_par},~\eqref{eps01} and~\eqref{aa1} -- \eqref{aa4} to obtain
\begin{align*}
a_{13}^2+a_{24}^2 & < 0.5\varepsilon^2, \\
a_{14}^2 + a_{23}^2 & \leq (a_{14}')^2 + (a_{23}')^2 < 0.5\varepsilon^2, \\
(a_{12}')^2 + (a_{34}')^2 & < 1.00251\varepsilon^2 0.5\varepsilon^2 < 0.5013\varepsilon^4, \\
(a_{12}'')^2+(a_{34}'')^2 & < (1+1.0052\varepsilon^2)0.5013\varepsilon^4 < 0.50634\varepsilon^4, \\
(a_{13}'')^2+(a_{24}'')^2 & < 1.00251\varepsilon^2 0.5013\varepsilon^4 < 0.5026\varepsilon^6.
\end{align*}
Now, using~\eqref{eps01}, we can find bounds for $S(A^{(6r+2l)})$, $1\leq l\leq 4$,
\begin{align*}
S^2(A^{(6r+2)}) & = (a_{13}^2+a_{24}^2)+(a_{14}^2+a_{23}^2) < 0.5\varepsilon^2+0.5\varepsilon^2 = \varepsilon^2, \\
S^2(A^{(6r+4)}) & = (a_{12}')^2+(a_{34}')^2+(a_{14}')^2+(a_{23}')^2 < 0.5013\varepsilon^4+0.5\varepsilon^2 < 0.52\varepsilon^2, \\
S^2(A^{(6r+6)}) & = (a_{12}'')^2+(a_{34}'')^2+(a_{13}'')^2+(a_{24}'')^2 < 0.5114\varepsilon^4 \\
S^2(A^{(6r+8)}) & = (a_{13}'')^2+(a_{24}'')^2 < 0.5026\varepsilon^6.
\end{align*}
Since $r\geq r_0$ is arbitrary, the last set of relations implies that all $S(A^{(6r_0+2l)})$, $l\geq 0$, are smaller than $\varepsilon$.
Thus, all elements of the sequence $S(A^{(2r)})$, $r\geq r_0$, except finitely many of them, are within $\varepsilon$-ball around $0$, which proves the theorem.

Further more, the latest relations indicate the ultimate cubic convergence per cycle, irrespective of the multiplicities of the eigenvalues.
\end{proof}

We end the paper by addressing less interesting cases which we depict as
\begin{align*}
A & = {\small\left[\begin{array}{r|rrr}
        \Box & \diamond & \diamond & \diamond \\ \hline
        \diamond & \Box & \cdot & \cdot \\
        \diamond & \cdot & \Box & \cdot \\
        \diamond & \cdot & \cdot & \Box \\
      \end{array}\right]}, \
J= {\small\left[\begin{array}{r|rrr}
        1 & 0 & 0 & 0 \\ \hline
        0 & -1 & 0 & 0 \\
        0 & 0 & -1 & 0 \\
        0 & 0 & 0 & -1 \\
      \end{array}\right]} \quad \text{and} \\
A & ={\small\left[\begin{array}{rrr|r}
        \Box & \cdot & \cdot & \diamond \\
         \cdot & \Box & \cdot & \diamond \\
        \cdot & \cdot  & \Box & \diamond \\ \hline
        \diamond & \diamond & \diamond & \Box \\
      \end{array}\right]}, \
J= {\small\left[\begin{array}{rrr|r}
        1 & 0 & 0 & 0 \\
        0 & -1 & 0 & 0 \\
        0 & 0 & -1 & 0 \\ \hline
        0 & 0 & 0 & -1 \\
      \end{array}\right]}.
\end{align*}
We can first consider the process on
{\scriptsize
$\displaystyle
\left[\begin{array}{rrr}
\Box & \cdot & \cdot \\
\cdot & \Box & \cdot \\
\cdot & \cdot & \Box
\end{array}\right]$}.
It is a cyclic Jacobi method with perturbations coming from
the hyperbolic transformations. Using the theory from \cite[Section~5]{BH_Jac2015}
one can prove that the off-norm of that $3$ by $3$ symmetric matric will tend to
zero. Namely, when $n=3$ every pivot ordering belongs to the class $\mathcal{C}_{sp}^{(3)}$.
The rest of the proof uses the fact that both, hyperbolic angles and the
corresponding pivot elements tend to zero.

\section*{Acknowledgements}
The authors are thankful to the anonymous reviewers for their suggestions and comments.

\end{document}